\documentclass[a4paper,11pt]{article}

\usepackage[utf8]{inputenc}
\usepackage[T1]{fontenc}
\usepackage{geometry}
\usepackage{amsmath}
\usepackage{amssymb}
\usepackage[all]{xy}
\usepackage{mathenv}
\usepackage[dvips]{color}
\usepackage{float}
\usepackage{psfrag}
\usepackage{pgf,tikz}
\usetikzlibrary{arrows}

\newtheorem{thm}{Theorem}

\newtheorem{dft}{Definition}
\newtheorem{prop}{Proposition}

\newtheorem{cor}{Corollary}
\newtheorem{rmk}{Remark}

\newtheorem{conjecture}{Conjecture}

\usepackage[french,english]{babel}
\usepackage{amssymb}

\def\C{\mathbb{C}}
\def\Cal{\mathcal{C}}
\def\I{\mathcal{I}}
\def\J{\mathcal{J}}

\def\H{\mathcal{H}}

\def\R{\mathbb{R}}

\def\N{\mathbb{N}}

\def\ve{\varepsilon}

\def\XXint#1#2#3{{\setbox0=\hbox{$#1{#2#3}{\displaystyle\int}$}\vcenter{\hbox{$#2#3$}}\kern-.5\wd0}}

\newcommand{\diff}{\mathop{}\mathopen{}\mathrm{d}}
\newcommand{\Sp}{\mathbb{S}}
\newcommand{\be}{\begin{enumerate}}
\newcommand{\ee}{\end{enumerate}}

\title{A necessary condition for lower semicontinuity of line energies}
\author{Pierre Bochard and Antonin Monteil \footnote{Laboratoire de Math\'ematiques, Universit\'e Paris-Sud, b\^at. 425, 91405 Orsay, France. 
Email: pierre.bochard@math.u-psud.fr, antonin.monteil@math.u-psud.fr} }

\begin{document}
\selectlanguage{english}
\maketitle

\begin{abstract}
We are interested in some energy functionals concentrated on the discontinuity lines of 
divergence-free 2D vector fields valued in the circle $\Sp^1$. This kind of energy has been
introduced first by P. Aviles and Y. Giga in \cite{MR924423}. They show in particular that, with the cubic cost function $f(t)=t^3$, 
this energy is lower semicontinuous. 
In this paper, we construct a counter-example which excludes the lower semicontinuity 
of line energies for cost functions of the form $t^p$ with $0<p<1$. 
We also show that, in this case, the viscosity solution corresponding to a certain convex domain is not a minimizer.

\end{abstract}

\section{Introduction}
\subsection{Line energies}

Let $\Omega$ be a Lipschitz domain in $\R^2$. We are interested in measurable vector fields $m \colon \Omega \to \R^2$ such that
\begin{equation}\label{eikonal}
|m|=1\text{ a.e. and } \nabla \cdot m=0\text{ on }\Omega,
\end{equation}
where the second equation holds in the distributional sense. In the following, we will assume that $m$ is of bounded variation so as to be able to define its {\it jump line}. So, we consider the set
\begin{equation*}
A(\Omega):=\left\{ m \in BV(\Omega,\R^2)\,:\, |m|=1 \text{ a.e. and }\nabla \cdot m=0 \text{ on } \Omega \right\}.
\end{equation*}
Vector fields $m\in A(\Omega)$ are related to solutions of the \textit{eikonal equation} in $\Omega$. Let define the set
\begin{equation*}
S(\Omega):=\{\varphi \in \text{Lip}(\Omega)\,:\, |\nabla\varphi|=1 \text{ a.e. and } \nabla\varphi\in BV(\Omega)\}.
\end{equation*}
Then, given $m\in A(\Omega)$, there exists a scalar function $\varphi\in S(\Omega)$ such that
\[
 m(x)=(\nabla\varphi(x))^\perp \text{ a.e.},
\]
where $(\nabla\varphi)^\perp=R\nabla\varphi$ stands for the image of $\nabla\varphi$ by the rotation $R$ of angle $\pi/2$ 
centered at the origin in $\R^2$. 
Moreover, a function $\varphi\in \operatorname{Lip}(\Omega)$ satisfying $(\nabla\varphi)^\perp=m$ a.e. is
unique up to a constant and is called stream function. We are now able to define line energies:\\

\begin{dft}\label{lineenergy}
Let $f \colon [0,2] \to [0,+\infty]$ be a measurable scalar function.
Let $m\in A(\Omega)\subset BV(\Omega, \R^2)$. Then, there exists a $\H^1$-rectifiable jump line $J(m)$ oriented by a unit 
normal vector $\nu_x$ such that $m$ has traces $m_{\pm}(x)\in \Sp^1$ on each side of $J(m)$ for $\H^1$ a.e. $x \in J(m)$ 
(see \cite{MR1857292} for more details). 

Then, the energy associated with the jump cost $f$ is denoted by $\I_f$ and defined for $m\in A(\Omega)$ as follows:
\begin{equation*}
\I_f(m)=\displaystyle\int_{J(m)}f(|m_+-m_-|) \diff \H^1(x).
\end{equation*}
\end{dft}
$f$ is called the jump cost. Note that the divergence constraint on $m\in A(\Omega)$ implies that for $a.e.\ x\in J(m)$, $m_\pm(x)\in \Sp^1$ and $\nu_x$
satisfy the following condition (see figure \ref{Nonsci}):
\begin{equation*}
m_+(x) \cdot \nu_x=m_-(x) \cdot \nu_x \ .
\end{equation*}
Then, in the orthogonal basis $(\nu_x,\nu_x^\perp)$, there exists some angle $\theta$ such that $m_\pm =(\cos\theta,\pm\sin\theta)$ 
and the jump size is defined as
\[
t=|m_+-m_-|=2|\sin\theta|.
\]
Similarly, $\I_f$ can be interpreted as a functional of the stream function on the set $S(\Omega)$: Writing $m=(\nabla\varphi)^\perp\in BV(\Omega,\R^2)$, then $\I_f(m)=\J_f(\varphi)$ where
\begin{equation}\label{jf}
\forall \varphi\in S(\Omega),\quad \J_f(\varphi)=\displaystyle\int_{J(\nabla\varphi)}f(|(\nabla\varphi)_+-(\nabla\varphi)_-|) \diff \H^1(x).
\end{equation}
An interesting question is to find the minimizing structures of $\I_f$ if it exists. Remark that for this problem to be relevant,
we have to consider a constraint on the boundary otherwise all constant functions are minimizers. A natural choice is to minimize 
$\I_f$ along all configurations $m$ belonging to the set
\begin{equation*}\label{constraint_on_m}
A_0(\Omega):=\left\{m \in A(\Omega)\;:\; m \cdot n=0 \text{ a.e. on } \partial \Omega \right\},
\end{equation*}
where $n$ is the exterior unit normal vector of $\partial \Omega$. In terms of the stream function $\varphi$, 
this is equivalent to consider the set
\begin{equation*}
S_0(\Omega):=\left\{ \varphi \in S(\Omega)\,:\, \varphi = 0\text{ on }\partial \Omega \right\}.
\end{equation*}

\subsection{Related models}

The first example of such energy is due to P. Aviles. and Y. Giga. In \cite{MR1669225}, they have conjectured that if $f(t)=\frac{1}{3}t^3$, 
then $\I_f$ is the $\Gamma$-limit of the following Ginzburg-Landau type energy functional
\begin{equation*}\label{AG}
AG_\ve (u)=
\begin{cases}
\displaystyle \int_\Omega \ve|\nabla u|^2+\frac{1}{\ve}(1-|u|^2)^2 &if $u \in H^1(\Omega,\R^2)$ and $\nabla \cdot u=0$, \\
+\infty &otherwise,
\end{cases}
\end{equation*}
where $\Omega$ is a bounded open set in $\R^2$ and $\ve>0$ is some parameter.

%
For the $\Gamma$-convergence of functionals $AG_\ve$ to $\I_f$ with $f(t)=\frac{1}{3} t^3$, only partial results are shown. 
In \cite{MR1669225}, the authors have been able to prove the $\Gamma$-liminf property for the $L^1$ convergence using the notion of 
entropies related to the problem \eqref{eikonal} (see also \cite{MR1752602}). 
The strong compactness of finite energy sequences
has been proved by Ambrosio, De Lellis and Mantegazza in \cite{MR1731470} and by
De Simone, Kohn, M\"uller and Otto in \cite{MR1854999} using a compensated 
compactness method based on a new notion of regular entropy on $\R^2$.\\


The second model we want to address comes from the Ginzburg-Landau theory in thin film micromagnetics for some asymptotical 
regime (see \cite{MR2763028}). Given a bounded domain $\Omega\subset \R^2$ and a magnetization $m=(m_1,m_2,m_3) \colon \Omega \to \Sp^2$, where $\Sp^2$ stands
for the unit sphere in $\R^3$, the energy of $m$ is defined as follows
\begin{equation}\label{fbetaeps}
E_\varepsilon(m)=
\begin{cases}
\displaystyle \ve \int_\Omega |\nabla m|^2 +\frac{1}{\ve}\int_\Omega \phi (m)& if $\nabla m'=0$ where $m'=(m_1,m_2)$,\\
+\infty & otherwise,
\end{cases}
\end{equation}
where $\ve$ is a small parameter called exchange length and $\phi : \Sp^2 \to \R$ is some smooth function called anisotropy function.
%
%

If $\phi (m) = |m_3|^{\alpha}$ with $0<\alpha\leq 4$, only one-dimensional structures are expected and it is easy to compute 
what should be the limiting energy of functionals $E_\ve$ by a $1D$-analysis. As for the Modica-Mortola model for phase 
transition (\cite{MR0445362}), $E_\ve$ is expected to $\Gamma$-converge to $c\,\I_f$ for some $c>0$ where $f(t)=\frac{t^p}{p}$, 
$p=1+\frac{\alpha}{2}$ is the primitive of $\sqrt{\phi}$ vanishing at $0$. The case $\phi (m)=|m_3|^2$ was studied by R. Ignat and
B. Merlet in \cite{MR2915327} in which a compactness result 
was proved and a lower bound was found. However, the $\Gamma$-liminf property in the definition of $\Gamma$-convergence was established 
only for limiting $1D$ configurations of the form $m(x)=\pm \nu^\perp$ for $\pm x \cdot \nu >0$ with $\nu\in \Sp^1$ (see figure \ref{Nonsci} 
with $\theta_0=\pi/2$).

\subsection{Lower semicontinuity, Viscosity solution}
As explained above, some of the line energies $\I_f$ are conjectured to be the $\Gamma$-limit of functionals coming from 
micromagnetics in the space $X=L^1$. If that is the case, $\I_f$ has to satisfy the following lower semicontinuity property:

\begin{dft}\label{lsc}
Let $F: X\to [0,+\infty]$ be a functional defined on some topological space $X$. $F$ is said to be lower semicontinuous or
l.s.c. if the following holds:
\[
\forall (x_n)_{n\geq 0}\subset X, \quad x_n \underset{n\to+\infty}{\longrightarrow}x \implies F(x)\leq \liminf_{n\to\infty}F(x_n).
\]
\end{dft}

Since this property strongly depends on the topology of the space $X$, we have to specify the choice we make for the study of 
line energies $\I_f$. \\
\indent First of all, due to the non convex constraint $|m|=1$, we need strong compactness in $L^1$. Moreover, since all 
the results of the previous part (compactness and $\Gamma$-liminf property) holds for the $L^1$ strong topology, it seems 
natural to consider the line energies $\I_f$ in the space $X=L^1$. \\

\indent However, since Definition \ref{lineenergy} uses the notion of trace of a function, another natural choice
would be $X=BV$ endowed with the weak topology which is a very common choice for phase transition problems. Unfortunately, in 
the general case, the space $BV$ is not adapted to our problem. Suppose $f(t)=t^p$ with $p>1$ for instance. Then finite energy configurations $m$ (i.e. $m_n \underset{n\to+\infty}{\longrightarrow} m$ 
in $L^1$ with $\I_f(m_n)\leq C<+\infty$) are not necessarily of bounded variation since the total variation of $m$ around 
its jump line cannot be controlled by $\int_{J(m)} |m_+-m_-|^p$ if $p>1$ (see \cite{MR1731470}). That is why we need a subspace of solutions of the problem \eqref{eikonal} included in $L^1(\Omega)$ 
(and containing $BV$) because of the non convex constraint $|m|=1$ such that we are still able to define a jump 
line $J(m)$ and traces $m_\pm$. This is done in \cite{MR1985613} 
where a regularity result is shown for solution of \eqref{eikonal} with bounded "entropy production".\\ 

Note that if $X$ and $Y$ are two topological spaces such that $Y$ is continuously embedded in
$X$ and $F: X\to [0,+\infty]$ is l.s.c. in $X$ then the restriction of $F$ to $Y$ is l.s.c. in $Y$. 
In this paper, we only want to prove a necessary condition for functionals $\I_f$ to be l.s.c. We then prefer 
to restrict our analysis to $BV$ functions (see remark \ref{rmklsc}). \\

In the case where $f(t)=t^p$ for some $p>0$, only partial results are known. In \cite{MR1731470}, the following is conjectured:

\begin{conjecture}
\label{conjecture}Let $\overline{\I_f}$ be the relaxation of $\I_f$ (only defined on the space $BV$) in $L^1$:
\begin{equation*}\label{Ifbar}
\overline{\I_f}(m)=\operatorname{Inf}
\left\{\liminf_{n \to +\infty} \I_f(m_n)\ :\ m_n\in BV\text{ and }m_n \underset{n\to+\infty}{\longrightarrow} m \text{ in }L^1 \right\}.
\end{equation*}
If $f(t)=t^p$ with $1\leq p\leq 3$ then $\overline{\I_f}$ is l.s.c. for the strong topology in $L^1$.
\end{conjecture}
For $p>3$, this conjecture is false (see \cite{MR1731470}).
The case $p=3$ has been studied by P. Aviles and Y. Giga in \cite{MR1669225}.  More recently 
the case $p=2$ has been proved by R. Ignat and B. Merlet in \cite{MR2915327}. They also proved that Conjecture \ref{conjecture} holds true for $1\leq p\leq 3$ if one restricts to configurations $m$ such that the jump size is always lower than $\sqrt{2}$. Here we
are interested in the open case $p<1$.\\
\indent We point out that line energies associated with the cost $f(t)=t^p$ with $1\leq p\leq 3$
correspond exactly to the expected $\Gamma$-limits of functionals \eqref{fbetaeps} when $\phi (m)=|m_3|^\alpha$ with
$0 <\alpha \leq 4$ where Bloch walls seem to be optimal. This is quite natural since when 2D 
structures, as cross tie wall or zigzag wall for instance, have less energy than Bloch walls, the $\Gamma$-limit 
of these functionals may be non lower semicontinuous. In the next part, we are going to give a 2D construction which 
gives some necessary condition on $f$ for $\I_f$ to be l.s.c. This condition excludes cost functions of the form $f(t)=t^p$ with $p<1$; the 
proof is based on a construction in the spirit of \cite{MR1932944} and \cite{MR1809740}.

\begin{thm}\label{CNLSC}
Let $f \colon [0,2] \to [0,+\infty ]$. Let $\Omega$ be an open and bounded non empty subset of $\R^2$. 
Assume that $\I_f$ is lower semicontinuous in $X=BV(\Omega,\Sp^1)$ endowed with the weak topology.
Then $f$ is lower semicontinuous and we have
\begin{equation}\label{CNF}
\limsup_{t\to 0}\ \frac{f(t)}{t} \leq 2\ \limsup_{t \to 2}\ f(t).
\end{equation}
\end{thm}

\begin{rmk}
The fact that the lower semicontinuity of $\mathcal{I}_f$ implies the lower semicontinuity of $f$ has already been proved in 
\cite{MR2915327}. The main new point here is the condition \eqref{CNF}. 
\end{rmk}
\begin{rmk}\label{rmklsc}
Theorem \ref{CNLSC} is stronger than an equivalent formulation in which $BV$ is replaced by some Banach space $X$ such that $BV$ is
continuously embedded in $X$ and where $\I_f$ is replaced by its relaxation in $X$.
\end{rmk}
As we will see, the lower semicontinuity of functionals $\I_f$ is closely related to the following question: Is the viscosity solution a 
minimizer of $\I_f$? In \cite{MR2915327}, the authors address the following conjecture
\begin{conjecture}
Assume that $\overline{\I_f}$ is l.s.c. in $L^1$ and that $\Omega$ is convex. Then $(\nabla\varphi_0)^\perp$ is a global
minimizer of $\I_f$ where 
$\varphi_0 (x)=\operatorname{dist}(x,\partial\Omega)$.
\end{conjecture}
For a regular domain $\Omega$ the distance function $\varphi_0(x)=\operatorname{dist}(x,\partial\Omega)$ belongs to
$S_0(\Omega)$ and $(\nabla\varphi_0)^\perp$ is the viscosity solution of the problem \eqref{eikonal}. In particular, 
if $\Omega$ is convex, $\varphi_0$ is concave and $-D^2 \varphi_0$ is a positive vectorial radon measure. In \cite{MR1731470}, 
the authors give a microstructure which shows that the viscosity solution is not a minimizer if $f(t)=t^p$ with $p>3$. As explained below, we are going 
to give a structure with lower energy than the viscosity solution for $p<1$.

\begin{prop}\label{CNV}
Let $f : [0,2] \to [0, +\infty]$. There exists a convex domain $\Omega$ such that the following holds. Let $\varphi_0\in S_0(\Omega)$ be the distance function $\varphi_0 (x)=\operatorname{dist}(x,\partial\Omega)$. Assume that
$\varphi_0$ is a minimizer of $\J_f$ defined by \eqref{jf}. Then f satisfies \eqref{CNF}.
\end{prop}

\begin{cor}\label{CorCNV}
There exists a convex domain $\Omega$ such that the viscosity solution is not a minimizer of $\I_f$ if $f(t)=t^p$ with $p\in [0,1[$.
\end{cor}

\section{Construction of a competitor of the viscosity solution}

In order to obtain the inequality \eqref{CNF}, we have to construct a domain $\Omega$ on which the jump size $t=|m_+-m_-|$ 
of the viscosity solution along its singular set is very small. Then, we find a competitor whose jump size $t$ is close 
to the maximal possible value $t=2$. In other words, we want to substitute small jumps by large ones.\\

We will use the polar coordinates $(r,\theta)$, $r\geq 0$, $\theta\in [- \pi, \pi]$ and we will identify 
$\R^2$ and $\C$ with the usual bijection. Let $D$ be the unit disk and $\Cal$ be its boundary. \\

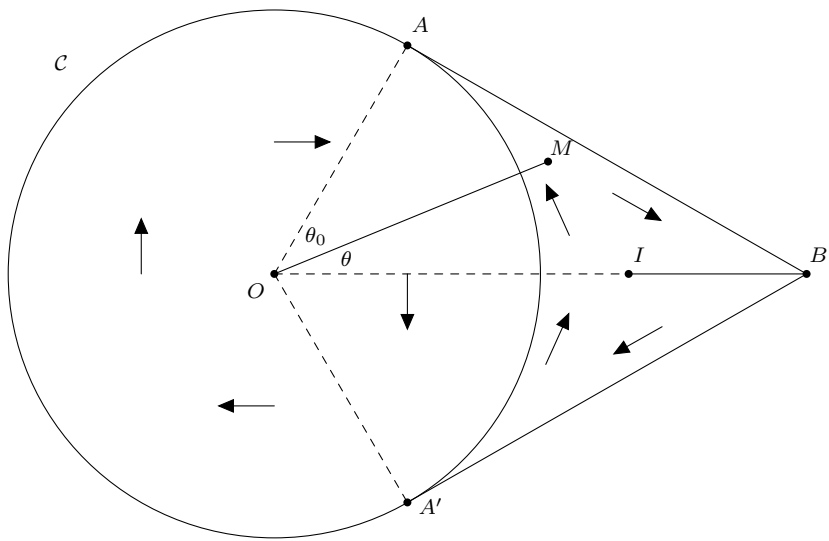
\begin{figure}[h]\center
\begin{tikzpicture}[line cap=round,line join=round,>=triangle 45,x=1.0cm,y=1.0cm]
\clip(-4,-4) rectangle (7.5,4);
\draw(0,0) circle (3.5cm);
\draw[dashed] (0,0)-- (1.75,3.03);
\draw [dashed] (0,0)-- (1.75,-3.03);
\draw[smooth,samples=100,domain=0.0:pi/3] plot[parametric] function{7*cos(t)/(1+cos(t-pi/3)),7*sin(t)/(1+cos(t-pi/3))};
\draw[smooth,samples=100,domain=0.0:pi/3] plot[parametric] function{7*cos(t)/(1+cos(t-pi/3)),-7*sin(t)/(1+cos(t-pi/3))};
\draw[smooth,samples=100,domain=0.0:pi/8] plot[parametric] function{0.8*cos(t),0.8*sin(t)};
\draw[smooth,samples=100,domain=0.0:pi/3] plot[parametric] function{0.5*cos(t),0.5*sin(t)};
\draw (14/3,0)-- (7,0);
\draw[dashed] (0,0)--(14/3,0);
\draw (1.75,3.03)-- (7,0);
\draw (1.75,-3.03)-- (7,0);
\draw [->] (-1.75,0) -- (-1.75,0.75);
\draw [->] (0,1.75) -- (0.75,1.75);
\draw [->] (1.75,0) -- (1.75,-0.75);
\draw [->] (0,-1.75) -- (-0.75,-1.75);
\draw [->] (4.45,1.07) -- (5.1,0.7);
\draw [->] (5.1,-0.7) -- (4.45,-1.07);
\draw (0,0)-- (3.6,1.49);
\draw [->] (3.88,0.51) -- (3.57,1.2);
\draw [->] (3.57,-1.2) -- (3.88,-0.51);
\begin{scriptsize}
\fill [color=black] (0,0) circle (1.5pt);
\draw[color=black] (-0.24,-0.22) node {$O$};
\fill [color=black] (1.75,3.03) circle (1.5pt);
\draw[color=black] (1.92,3.3) node {$A$};
\fill [color=black] (1.75,-3.03) circle (1.5pt);
\draw[color=black] (2.06,-3.08) node {$A'$};
\fill [color=black] (7,0) circle (1.5pt);
\draw[color=black] (7.16,0.26) node {$B$};
\fill [color=black] (4.66,0) circle (1.5pt);
\draw[color=black] (4.8,0.26) node {$I$};
\fill [color=black] (3.6,1.49) circle (1.5pt);
\draw[color=black] (3.78,1.67) node {$M$};
\draw[color=black] (0.95,0.2) node{$\theta$};
\draw[color=black] (0.55,0.5) node{$\theta_0$};
\draw[color=black] (-2.8,2.8) node{$\mathcal{C}$};
\end{scriptsize}
\end{tikzpicture}
\caption{\label{microstruct}The domain $\Omega$ and the microstructure $m$} 
\end{figure}

Let  $\theta_0$ be a fixed angle in $]0,\pi/2[$  and define the two points $A=e^{i\theta_0}$ 
and $A'=e^{-i\theta_0}$ on the circle $\Cal$. Define also $T_A$ (resp. $T_{A'}$) the tangent 
to the circle $\Cal$ at the point $A$ (resp. $A'$). We consider the domain $\Omega$ delimited by the large arc $\{e^{i\theta}\;:\;|\theta|>\theta_0\}$, 
$T_A$ and $T_{A'}$ (see figure \ref{microstruct}). In other words $\Omega$ is the interior of the convex envelope of 
$\Cal \cup \{B\}$ where $B=T_A\cap T_{A'}$.
Define also 
$\Omega_0=\Omega\cap\{|\theta|<\theta_0\text{ and }r>0\}$ and  $\Gamma=\partial\Omega\cap\partial\Omega_0=[AB]\cup[A'B]$.\\

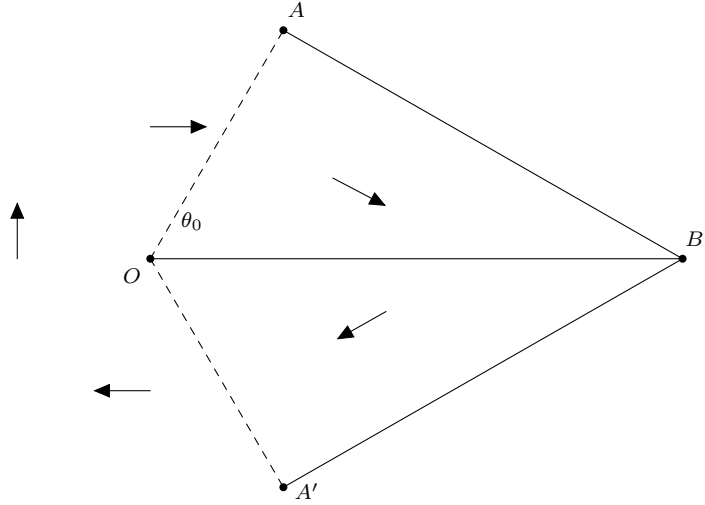
\begin{figure}[h!]
\centering
\begin{tikzpicture}[line cap=round,line join=round,>=triangle 45,x=1.0cm,y=1.0cm]
\clip(-4,-4) rectangle (7.5,4);
\draw[dashed] (0,0)-- (1.75,3.03);
\draw[dashed] (0,0)-- (1.75,-3.03);
\draw[smooth,samples=100,domain=0.0:pi/3] plot[parametric] function{0.5*cos(t),0.5*sin(t)};
\draw[smooth,samples=100,domain=pi/3:(5*pi/3)] plot[parametric] function{3.5*cos(t),3.5*sin(t)};
\draw (0,0)-- (7,0);
\draw (1.75,3.03)-- (7,0);
\draw (1.75,-3.03)-- (7,0);
\draw [->] (-1.75,0) -- (-1.75,0.75);
\draw [->] (0,1.75) -- (0.75,1.75);
\draw [->] (0,-1.75) -- (-0.75,-1.75);
      \draw [->] (4.4-2,1.07) -- (5.1-2,0.7);
      \draw [->] (5.1-2,-0.7) -- (4.45-2,-1.07);
\begin{scriptsize}
\fill [color=black] (0,0) circle (1.5pt);
\draw[color=black] (-0.24,-0.22) node {$O$};
\fill [color=black] (1.75,3.03) circle (1.5pt);
\draw[color=black] (1.92,3.3) node {$A$};
\fill [color=black] (1.75,-3.03) circle (1.5pt);
\draw[color=black] (2.06,-3.08) node {$A'$};
\fill [color=black] (7,0) circle (1.5pt);
\draw[color=black] (7.16,0.26) node {$B$};
\draw[color=black] (0.55,0.5) node{$\theta_0$};
\end{scriptsize}
\end{tikzpicture}
\caption{\label{Nonsci}Viscosity solution $m_0$ on $\Omega$} 
\end{figure}

We now consider two solutions $\varphi_0$ and $\varphi$ in $S_0(\Omega)$ of the eikonal equation vanishing on the boundary:
\begin{itemize}
\item $\varphi_0$ is the usual distance function: $\forall x\in\Omega,\; \varphi_0(x)=\operatorname{dist}(x,\partial\Omega)$.
\item $\varphi$ defined by: $\forall x\in\Omega,\; \varphi (x)=\operatorname{dist}(x,\partial\Omega\cup\Cal)$.
\end{itemize}

We also denote by $m_0=(\nabla \varphi_0)^\perp$ and $m=(\nabla\varphi)^\perp$ the corresponding solutions of \eqref{eikonal}. 
Then $m_0, m\in A_0(\Omega)$.

We now compute $\I_f(m_0)$ and $\I_f(m)$ in order to prove that the function $\varphi$ has lower energy than $\varphi_0$ if $f(t)=t^p$
with $p<1$. 

\paragraph{Heuristic:}

The idea is that a small jump along a fixed length is replaced by big jumps on a small length : 
This will reduce the energy for subadditive power costs (i.e. $f(t)=t^p$ with $p<1$) which favor "small jumps". 
Let us give more details.\\\indent
For a small angle $\theta_0>0$, $m_0$ only presents small jumps:
$m_0$ is $\mathcal{C}^1$ out of segment $[OB]$ on which the jump size is $|m_0^+-m_0^-|=:t_0=2\sin(\theta_0)$. \\
\indent On the contrary, $m$ only presents "big" jumps: i.e. jumps whose size is close to $2$. The singular set of $m$ consists 
in $3$ different lines : $[IB]$ whose length is equivalent to $\theta_0^2$ and the two curves $\Cal\setminus\Cal_{\theta_0}$ and
$\gamma_{\theta_0}$ (defined below) on which the jump size tends to $2$ and the length of these lines is equivalent to $2\theta_0$.\\\indent
 As a result, the energy of $m_0$ is close to $f(2\sin\theta_0)$ while the energy of $m$ is close to $4\theta_0\times f(2)$. 
A necessary condition for $m_0$ to minimize $\I_f$ is then (see Proposition \ref{CNV})
$$\underset{t\to 0}{\limsup}f(t)/t\leq 2 f(2).$$ 
This excludes subadditive power costs. In the sequel, we are going to make precise computations so as to get more informations about the critical angle $\theta_0$.

\paragraph{Energy of $m_0$:}

The jump line of $m_0$ is the segment $[OB]$ and the traces of $m_0$ on each side of this line are given by
$m_{0,\pm}=-e^{i(\pi/2\pm\theta_0)}$. In particular,
\begin{equation*}\label{E(m_0)}
\I_f(m_0)=f(2\sin\theta_0)|OB|=\frac{f(2\sin\theta_0)}{\cos\theta_0}.
\end{equation*}
\paragraph{Energy of $m$:}

The jump line of $m$ is the union of the 3 curves:

\begin{itemize}
\item\label{c} $\Cal_{\theta_0}=\{e^{i\theta} \,:\,|\theta|<\theta_0\}$.
\item $\gamma_{\theta_0}:=\{z\in\Omega_0\;:\;d(z,\Cal_{\theta_0})=d(z,\Gamma)\}=\{z=re^{i\theta}\;:\;|\theta|<\theta_0,\, d(z,\Cal)
=d(z,\partial\Omega)\}.$
\item The segment $[IB]$ where $I=\gamma_{\theta_0}\cap[OB]$.

\end{itemize}
  
First, let us find a polar equation for the curve $\gamma_{\theta_0}$: Given $z=re^{i\theta}$ such that
$|\theta|<\theta_0$ and $r>1$ we have $d(z,\Cal_{\theta_0})=r-1$, it remains to compute $\lambda:=d(z,\Gamma)$. \\

Since $\Omega$ is symmetric with respect to the axis $(OB)$, one can restrict to the case 
$M=r\, e^{i\theta}$ with $0<\theta<\theta_0$. So $\lambda :=d(z,\Gamma)=|z-P|$ where $P$ is the orthogonal projection of
$M=re^{i\theta}$ on the segment $[AB]$ : $P$ should satisfy $\overrightarrow{MP}=\lambda\, 
\overrightarrow{OA}=\lambda\, e^{i\theta_0}$ and $\overrightarrow{MP} \cdot \overrightarrow{AP}=0$. We then compute
\begin{align*}
\overrightarrow{MP} \cdot \overrightarrow{AP} = &  \overrightarrow{MP} \cdot [\overrightarrow{AO}+\overrightarrow{OM}+\overrightarrow{MP}], \\
 = & \Re \{\lambda\,e^{-i\theta_0}\, (-e^{i\theta_0}+r\,e^{i\theta}+\lambda\,e^{i\theta_0})\}, \\
 = & \lambda [-1+r\cos(\theta_0-\theta)+\lambda].
\end{align*}
Since $\overrightarrow{MP} \cdot \overrightarrow{AP}=0$, this implies $\lambda=MP=1-r\cos(\theta_0-\theta)$. Then we have $z\in\gamma_{\theta_0}$ if and only if $r-1=1-r\cos(\theta_0-\theta)$ and the polar equation of the curve $\gamma_{\theta_0}$ is given by
\begin{equation*}
r(\theta)=\displaystyle\frac{2}{1+\cos(\theta_0-|\theta|)}\ ;\ -\theta_0<\theta<\theta_0\quad .
\end{equation*}
Now, we can compute the energy of $m$ along the curve $\gamma_{\theta_0}$:
\begin{itemize}
\item $\diff \gamma(\theta)=\sqrt{r(\theta)^2+r'(\theta)^2} \diff \theta$ where we find $r'(\theta)=
\displaystyle\frac{-2\sin(\theta_0-\theta)}{(1+\cos(\theta_0-\theta))^2}$. Introducing the notation $\alpha=\theta_0-\theta$, we obtain
\[
\diff \gamma(\theta)=2\displaystyle\frac{\sqrt{(1+\cos\alpha)^2+\sin^2\alpha}}{(1+\cos\alpha)^2} \diff \theta
=2\displaystyle\frac{\sqrt{2(1+\cos\alpha)}}{(1+\cos\alpha)^2} \diff \theta=\displaystyle\frac{4\cos(\alpha/2) }{(2\cos^2(\alpha/2))^2} \diff \theta .
\]
So $\diff \gamma$ reads
\[
\diff \gamma(\theta)=\cos^{-3}(\alpha/2)\diff \theta .
\]

\item The size of the jump at the point $\gamma(\theta)$ is given by
\[
t(\theta)=|m_+-m_-|=|e^{i(\theta_0+\pi/2)}+e^{i(\theta+\pi/2)}|=|e^{i(\theta_0-\theta)}+1|.
\]

Using once again the notation $\alpha=\theta_0-\theta$, this gives
\[
t(\theta)=\sqrt{(\cos\alpha+1)^2+\sin^2\alpha}=\sqrt{2(1+\cos\alpha)}=2\cos(\alpha/2).
\]

\item We conclude that the energy of $m$ induced by the jump line $\gamma_{\theta_0}$ is given by
\begin{equation*}
\I_f^1(m) =\displaystyle\int_{-\theta_0}^{\theta_0}\displaystyle\frac{f[2\cos(\alpha/2)]}{\cos^3(\alpha/2)} \diff \alpha .
\end{equation*}
\end{itemize}
The energy concentrated on the arc $\Cal_{\theta_0}$ is
\begin{equation*}
\I_f^2(m)=f(2)\,\H^1(\Cal_{\theta_0})=2\theta_0 \,f(2).
\end{equation*}
Finally, we compute the energy on the line $[IB]$:
\begin{equation*}
\I_f^3(m)=f(2\sin\theta_0)|IB|.
\end{equation*}
If the distance function is a minimizer of $\I_f$ we should have
\[
\I_f(m)-\I_f(m_0)\geq 0.
\]
Now, the preceding equations yields
\begin{align*}
\I_f(m)-\I_f(m_0) =& \I_f^1(m)+\I_f^2(m)+\I_f^3(m)-\I_f(m_0),\\
= & \int_{-\theta_0}^{\theta_0}\frac{f[2\cos(\alpha/2)]}{\cos^3(\alpha/2)} \diff \alpha +2\theta_0\, f(2)+\left(|IB|-|OB| \right) f(2\sin\theta_0).
\end{align*}
Since $|IB|-|OB|=-|OI|=-r(0)=-\displaystyle \frac{1}{\cos^2(\theta_0/2)}$, this gives
\begin{equation*}\label{diff}
\I_f(m)-\I_f(m_0)=\int_{-\theta_0}^{\theta_0}\displaystyle\frac{f[2\cos(\alpha/2)]}{\cos^3(\alpha/2)} \diff \alpha
+2\theta_0 f(2)-\displaystyle\frac{f(2\sin\theta_0)}{\cos^2(\theta_0/2)}.
\end{equation*}
Hence, if $m_0$ is a minimizer of $\I_f$, the following condition should be satisfied:
\begin{align*}
\displaystyle\frac{f(2\sin\theta_0)}{2\sin\theta_0} &\leq\ \displaystyle\frac{\theta_0\cos^2(\theta_0/2)}{\sin\theta_0}
\left[ \frac{1}{\theta_0}\displaystyle\int_{0}^{\theta_0}\displaystyle\frac{f[2\cos(\alpha/2)]}{\cos^3(\alpha/2)}\diff \alpha+ f(2)\right] ,\\
 &\leq\  \displaystyle\frac{\theta_0}{\sin\theta_0\cos(\theta_0/2)}\times 2\ \sup\{f(t) \,:\,  2 \cos(\theta_0/2)\leq t\leq 2\}.
\end{align*}
Finally, taking the $\limsup$ for $\theta_0\to 0$ in the preceding equation leads to \eqref{CNF}:
\begin{equation*}
\limsup_{t \to 0} \frac{f(t)}{t} \leq 2 \limsup_{t \to 2} f(t).
\end{equation*}
This proves Proposition \ref{CNV} and corollary \ref{CorCNV} follows from the fact that the preceding inequality holds 
false for $f(t)=t^p$ with $p<1$. Note that in this case, we get something more precise than Proposition \ref{CNV}:

\begin{prop}
There exists $\theta_0\in ]0,\pi/2[$ only depending on $p$ such that for all\\ $\theta\in ]-\theta_0,\theta_0[$, the viscosity 
solution is not a minimizer of $\I_f$ on $\Omega_{\theta}$ where $\Omega_\theta$ is the convex set constructed in the previous part 
($\theta$ being the angle $(\overrightarrow{OB},\overrightarrow{OA})$).
\end{prop}

\section{Lower semicontinuity of line energies, proof of Theorem \ref{CNLSC}.}

The fact that if $\mathcal{I}_f$ is l.s.c. then $f$ is l.s.c. can be found in \cite{MR2915327} (Proposition 1).
In this section we prove that \eqref{CNF} is a necessary condition for $\I_f$ to be lower semicontinuous with respect to the weak
convergence in $BV$ on bounded open subsets of $\R^2$.

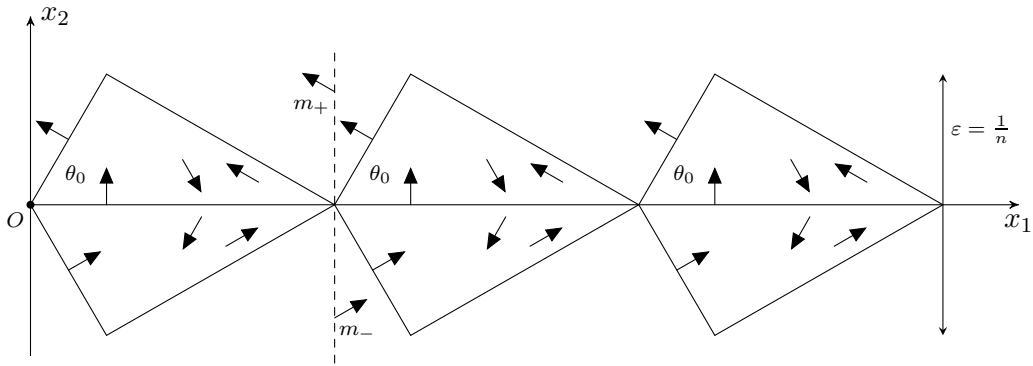
\begin{figure}[h]
\begin{tikzpicture}[line cap=round,line join=round,>=triangle 45,x=1.0cm,y=1.0cm]
\clip(-4.5,-3) rectangle (10,3);

\draw[smooth,samples=100,domain=0:pi/3] plot[parametric] function{-4+4*cos(t)/(1+cos(t-pi/3)),-4*sin(t)/(1+cos(t-pi/3))};
\draw[smooth,samples=100,domain=0:pi/3] plot[parametric] function{-4+4*cos(t)/(1+cos(t-pi/3)),4*sin(t)/(1+cos(t-pi/3))};
\draw (-4+1,2*0.866025)-- (-4+4,0);
\draw (-4+1,-2*0.866025)-- (-4+4,0);
\draw (-4,0)-- (-4+1,2*0.866025);
\draw (-4,0)-- (-4+1,-2*0.866025);
\draw[smooth,samples=100,domain=-pi/3:pi/3] plot[parametric] function{-4+2*cos(t),2*sin(t)};
\draw[smooth,samples=100,domain=0:pi/3] plot[parametric] function{-4+0.5*cos(t),0.5*sin(t)};
\draw [->] (-4+0.5,0.866025) -- (-4+0.5-1/2*0.866025,0.866025+1/2*0.5);
\draw [->] (-4+1,0) -- (-4+1,0.5);
\draw [->] (-4+0.5,-0.866025) -- (-4+0.5+1/2*0.866025,-0.866025+1/4);
\draw [->] (-4+3,0.3) -- (-4+2.57,0.55);
\draw [->] (-4+2.57,-0.55) -- (-4+3,-0.3);
\draw [->] (-4+2,0.6) -- (-4+2.25,0.16);
\draw [->] (-4+2.25,-0.16) -- (-4+2,-0.6);

\draw[smooth,samples=100,domain=0:pi/3] plot[parametric] function{4*cos(t)/(1+cos(t-pi/3)),-4*sin(t)/(1+cos(t-pi/3))};
\draw[smooth,samples=100,domain=0:pi/3] plot[parametric] function{4*cos(t)/(1+cos(t-pi/3)),4*sin(t)/(1+cos(t-pi/3))};
\draw (1,2*0.866025)-- (4,0);
\draw (1,-2*0.866025)-- (4,0);
\draw (0,0)-- (1,2*0.866025);
\draw (0,0)-- (1,-2*0.866025);
\draw[smooth,samples=100,domain=-pi/3:pi/3] plot[parametric] function{2*cos(t),2*sin(t)};
\draw[smooth,samples=100,domain=0:pi/3] plot[parametric] function{0.5*cos(t),0.5*sin(t)};
\draw [->] (0.5,0.866025) -- (0.5-1/2*0.866025,0.866025+1/2*0.5);
\draw [->] (1,0) -- (1,0.5);
\draw [->] (0.5,-0.866025) -- (0.5+1/2*0.866025,-0.866025+1/4);
\draw [->] (3,0.3) -- (2.57,0.55);
\draw [->] (2.57,-0.55) -- (3,-0.3);
\draw [->] (2,0.6) -- (2.25,0.16);
\draw [->] (2.25,-0.16) -- (2,-0.6);

\draw[smooth,samples=100,domain=0:pi/3] plot[parametric] function{4+4*cos(t)/(1+cos(t-pi/3)),-4*sin(t)/(1+cos(t-pi/3))};
\draw[smooth,samples=100,domain=0:pi/3] plot[parametric] function{4+4*cos(t)/(1+cos(t-pi/3)),4*sin(t)/(1+cos(t-pi/3))};
\draw (4+1,2*0.866025)-- (4+4,0);
\draw (4+1,-2*0.866025)-- (4+4,0);
\draw (4,0)-- (4+1,2*0.866025);
\draw (4,0)-- (4+1,-2*0.866025);
\draw[smooth,samples=100,domain=-pi/3:pi/3] plot[parametric] function{4+2*cos(t),2*sin(t)};
\draw[smooth,samples=100,domain=0:pi/3] plot[parametric] function{4+0.5*cos(t),0.5*sin(t)};
\draw [->] (4+0.5,0.866025) -- (4+0.5-1/2*0.866025,0.866025+1/2*0.5);
\draw [->] (4+1,0) -- (4+1,0.5);
\draw [->] (4+0.5,-0.866025) -- (4+0.5+1/2*0.866025,-0.866025+1/4);
\draw [->] (4+3,0.3) -- (4+2.57,0.55);
\draw [->] (4+2.57,-0.55) -- (4+3,-0.3);
\draw [->] (4+2,0.6) -- (4+2.25,0.16);
\draw [->] (4+2.25,-0.16) -- (4+2,-0.6);

\draw[dashed] (0,-2.1)--(0,2.1);
\draw [->] (0,1.5)--(0-0.5*0.866025,1.5+0.5*1/2);
\draw [->] (0,-1.5)--(0+0.5*0.866025,-1.5+0.5*1/2);

\draw [->,>=stealth'] (-4,-2) -- (-4,2.5) node[right] {$x_2$};
\draw [->,>=stealth'] (-4,0) -- (9,0) node[below] {$x_1$};

\draw [<->,>=stealth] (8,-2*0.866025)--(8,2*0.866025);

\begin{scriptsize}
\fill [color=black] (-4,0) circle (1.5pt);
\draw[color=black] (-4.2,-0.2) node {$O$};
\draw[color=black] (0.6-4,0.4) node {$\theta_0$};
\draw[color=black] (0.6,0.4) node {$\theta_0$};
\draw[color=black] (0.6+4,0.4) node {$\theta_0$};
\draw[color=black] (-0.3,1.3) node{$m_+$};
\draw[color=black] (0.3,-1.7) node{$m_-$};
\draw[color=black] (8.5,1) node{$\varepsilon= \frac{1}{n}$};
\end{scriptsize}
\end{tikzpicture}
\caption{\label{align} The microstructure $m_n$} 
\end{figure}

The key is to use the construction $m\in A(\Omega)$ depending on $\theta_0$ of the first part by restriction to $\Omega_0$ 
(See figure \ref{Nonsci}.).  The $1D$ transition defined by \eqref{1Dtransition} corresponds to the viscosity solution $m_0$ 
of the previous part. Given a small parameter $\epsilon >0$, it will costs less energy to substitute the $1D$ transition around 
its jump line by the microstructure $m$ rescaled at the level $\epsilon$ (see figure \ref{align}). \\

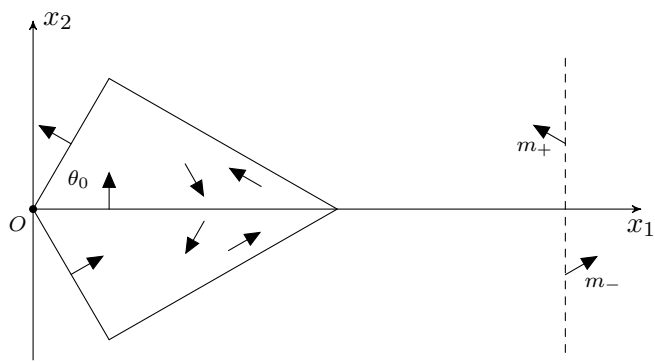
\begin{figure}[H]
\centering
\begin{tikzpicture}[line cap=round,line join=round,>=triangle 45,x=1.0cm,y=1.0cm]
\clip(-1,-3) rectangle (9,3);

\draw[smooth,samples=100,domain=0:pi/3] plot[parametric] function{4*cos(t)/(1+cos(t-pi/3)),-4*sin(t)/(1+cos(t-pi/3))};
\draw[smooth,samples=100,domain=0:pi/3] plot[parametric] function{4*cos(t)/(1+cos(t-pi/3)),4*sin(t)/(1+cos(t-pi/3))};
\draw (1,2*0.866025)-- (4,0);
\draw (1,-2*0.866025)-- (4,0);
\draw (0,0)-- (1,2*0.866025);
\draw (0,0)-- (1,-2*0.866025);
\draw[smooth,samples=100,domain=-pi/3:pi/3] plot[parametric] function{2*cos(t),2*sin(t)};
\draw[smooth,samples=100,domain=0:pi/3] plot[parametric] function{0.5*cos(t),0.5*sin(t)};
\draw [->] (0.5,0.866025) -- (0.5-1/2*0.866025,0.866025+1/2*0.5);
\draw [->] (1,0) -- (1,0.5);
\draw [->] (0.5,-0.866025) -- (0.5+1/2*0.866025,-0.866025+1/4);

\draw [dashed] (7,2)--(7,-2);
\draw [->] (7,0.87) -- (6.57,1.12);
\draw [->] (7,-0.87) -- (7.43,-0.62);
\draw [->] (3,0.3) -- (2.57,0.55);
\draw [->] (2.57,-0.55) -- (3,-0.3);
\draw [->] (2,0.6) -- (2.25,0.16);
\draw [->] (2.25,-0.16) -- (2,-0.6);
\draw [->,>=stealth'] (0,-2) -- (0,2.5) node[right] {$x_2$};
\draw [->,>=stealth'] (0,0) -- (8,0) node[below] {$x_1$};

\begin{scriptsize}
\fill [color=black] (0,0) circle (1.5pt);
\draw[color=black] (-0.2,-0.2) node {$O$};
\draw[color=black] (0.6,0.4) node {$\theta_0$};
\draw[color=black] (6.6,0.8) node{$m_{+}$};
\draw[color=black] (7.5,-1) node{$m_{-}$};
\end{scriptsize}
\end{tikzpicture}
\caption{\label{Nonsci}The vector field $m$ on the left and the 1D-transition $m_0$ on the right} 
\end{figure}

We are going to prove Theorem \ref{CNLSC} when $\Omega=(0,1)\times (-1,1)$. The general case will follow easily. Fix $\theta_0\in ]0,\pi/2[$ and define the $1D$ transition $m_0$ for a.e. $x_1\in (0,1)$ and $x_2\in \R$ by
\begin{equation}\label{1Dtransition}
m_0(x_1,x_2)=m_\pm:=(\mp\sin\theta_0, \cos\theta_0)\text{   if   }\pm x_2>0.
\end{equation}
Then, let us consider the vector field $m=m_{\theta_0}$ of the preceding section restricted to $\Omega_0$ and define the rescaled and extended 
vector field $\tilde{m}$ for $x_1\in (0,1)$ and $x_2\in \R$: 
\begin{equation*}
 \tilde{m}(x_1,x_2)=
 \left\{ \begin{array}{ll}
 -m\left((\cos\theta_0)^{-1}\; x_1,(\cos\theta_0)^{-1}\; x_2\right) & \text{if} 
 \left((\cos\theta_0)^{-1}\,x_1,(\cos\theta_0)^{-1}\,x_2\right)\in\Omega_0, \\ 
 m_0(x_1,x_2) & \text{otherwise.} 
 \end{array} \right.
\end{equation*}
Note that $\tilde{m}$ belongs to $A(\Omega)$ and is continuous up to the boundary, $\tilde{m}\in\mathcal{C}\left(\overline{\Omega}_0\right)$.
Then, let $n$ be a positive integer and define $m_n\in A(\Omega)$ by aligning n times the vector field $\tilde{m}$ (see figure \ref{align}).
More precisely, for $0\leq i<n$ and $x=(x_1,x_2)\in \Omega$ such that $i/n\leq x_1 <(i+1)/n$, define
\[
m_n(x_1,x_2)=\tilde{m}(n\, x_1-i,n\, x_2).
\]
We have $m_n(x_1,x_2)=m_0(x_1,x_2)$ for $|x_2|>1/n$ and for all $x\in \Omega,\, |m_n(x)|=1$. Consequently, $(m_n)_{n>0}$ converge to 
$m_0$ in $L^1(\Omega)$.
Moreover, $|m_n|_{BV(\Omega)}=|\tilde{m}|_{BV(\Omega)}$ so that $(m_n)_{n>0}$ is bounded
in $BV(\Omega)$ and weakly converges to $m_0$.

Since $m_n$ is obtained by scaling a fixed structure, it is easy to see that $\I_f(m_n)$ is constant. Indeed, 
$\I_f(m_n)=n\times 1/n\ \I_f(\tilde{m})=\I_f(\tilde{m})$. That is why we obtain the following condition: 
assuming $\I_f$ is $l.s.c.$,
\[
\I_f(m_0)=f(2\,\sin\theta_0)\leq\underset{n\to\infty}{\liminf}\ \I_f(m_n)=\I_f(\tilde{m}).
\]
In other words, the viscosity solution costs less energy than the construction $m_{\theta_0}$ of the preceding part. For this
reason, we obtain exactly the same necessary condition \eqref{CNF} and this ends the proof when $\Omega=(0,1)\times (-1,1)$. 

In the general case, let $D$ be an horizontal line such that $D \cap \Omega \neq \emptyset$. Up to a translation of $\Omega$, we can assume that $0\in D$, i.e. $D=\{x\in\R^2\;:\; x_2=0\}$. Fix $\theta_0\in ]0,\pi/2[$ and define $\tilde{m}_0\in\operatorname{BV}(\Omega,\Sp^1)$ by:
\begin{equation}\label{mtildeomega}
\text{for a.e. }x\in\Omega,\; \tilde{m}_0(x_1,x_2):=m_\pm=(\mp\sin\theta_0, \cos\theta_0)\quad\text{if}\quad\pm x_2>0\quad .
\end{equation}
Let $q \in \R$ and $r>0$ and define $\Omega_{q,r}=(q,q+r) \times (-r,r)$. $n$ being fixed,
we obtain a microstructure $m_n^{q,r}$ defined on $\Omega_{q,r}$
by rescaling and translating our microstructure $m_n$ defined on $(0,1) \times (-1,1)$. More precisely,
for a.e. $x\in\Omega_{q,r}$,
\[
 m_n^{q,r}(x_1,x_2):= m_n\left( \frac{x_1-q}{r},\frac{x_2}{r} \right). 
\]
Now, there exists two sequences $(q_i)_{i \in \N},\, (r_i)_{i\in \N}$ such that
\[
 D \cap \Omega \subset \bigcup_{i=0}^{\infty}\overline{\Omega_{q_i,r_i}}\quad ,
\]
where for $i\neq j$, $ \Omega_{q_i,r_i} \cap \Omega_{q_j,r_j} =\emptyset$.
We can now define the microstructure $\tilde{m}_n$ for a.e. $x$ belonging to the whole domain $\Omega$ by
\begin{equation*}
 \tilde{m}_n(x)=
\begin{cases}
m_n^{q_i,r_i}(x) &if $x \in \Omega_{q_i,r_i}$ for some $i\in\N$, \\
\tilde{m}_0(x) &otherwise,
\end{cases}
\end{equation*}
where $\tilde{m}_0$ has been defined in \eqref{mtildeomega}. It is clear that $\tilde{m}_n \xrightarrow[n \to \infty]{} \tilde{m}_0$ in $L^1(\Omega)$ and in the weak $BV$ sense. Furthemore,
\[
\I_f(\tilde{m}_0)=\operatorname{length}(\Omega\cap D)\ \I_f(m_0)\leq\underset{n\to\infty}{\liminf}\ \I_f(\tilde{m}_n)=\operatorname{length}(\Omega\cap D)\ \underset{n\to\infty}{\liminf}\ \I_f(m_n).
\]
Since $\Omega\cap D\neq \emptyset$ and $\Omega$ is bounded, one has $+\infty>\operatorname{length}(\Omega\cap D)>0$ which imply the thesis.
\paragraph{Aknowledgment.}
We want to thank our advisor Radu Ignat for pointing us this interesting problem and helpful discussions.

\newpage

\bibliographystyle{plain}

\end{document}